\documentclass{article}
\usepackage{amsmath}
\usepackage{amsfonts}
\usepackage{amsthm}
\newcommand{\RR}{\mathbb{R}}
\newtheorem{theorem}{Theorem}
\newtheorem*{thm1}{Theorem 1}
\newtheorem{lemma}{Lemma}
\newtheorem{cor}{Corollary}
\theoremstyle{definition}
\newtheorem{defn}{Definition}
\theoremstyle{remark}
\newtheorem*{acn}{Acnowledgements}

\DeclareMathOperator{\inv}{inv}
\DeclareMathOperator{\sgn}{sgn}
\DeclareMathOperator{\pat}{pat}
\DeclareMathOperator{\lat}{lat}
\begin{document}
\title{Symplectic tensor invariants, wave graphs and S-tris 
\thanks{This research was partially supported by 
ONR grant N00014-97-1-0505.}}

\author{Aleksandrs Mihailovs\\
Department of Mathematics\\
University of Pennsylvania\\
Philadelphia, PA 19104-6395\\
mihailov@math.upenn.edu\\
http://www.math.upenn.edu/$\sim$mihailov/
}
\date{\today}
\maketitle

\begin{abstract}
The spaces of invariants of tensor powers of the defining representation of 
$Sp(2n)$ are provided with the bases parametrized by symplectic wave graphs 
introduced here especially for this purpose. The proof utilizes a game similar to 
Tetris,  named here S-tris. This work continues my previous work \cite{M3} 
on the tensor invariants of $SL(n)$, wave graphs and L-tris.
\end{abstract}    
\setlength{\baselineskip}{1.5\baselineskip}

\section{Introduction}

Rumer, Teller and Weyl \cite{W} parametrized a basis of 
the subspace of $SL(2)$-invariants of 
$V^{\otimes m}=V \otimes \dots \otimes V$ ($m$ times), 
where $V$ is the two-dimensional linear space 
with the standard action of $SL(2)$, by  
$1$-regular outerplanar graphs, i.\ e.\ graphs with the vertices $1,2, \dots, m$, 
edges of which can be drawn in the upper half-plane without intersections.
They used slightly different graphs, drawn as a set of non-intersecting 
chords inside a disk, but after a conformal mapping of a disk onto the 
upper half plane and the indexing of the vertices in the increasing order, 
one gets the graphs described above. 

This theory was developed and 
applied to the percolation theory by Temperley and Lieb \cite{TL},
to the knots theory and invariants of 3-manifolds by Jones \cite{J}, 
Kauffman \cite{Kauf}, Kauffman and Lins \cite{KL}, Wenzl \cite{Wen}, 
Jaeger \cite{Ja}, Lickorish \cite{L}, Masbaum and Vogel \cite{MV} and others, 
to the quantum theory by Penrose \cite{P} and Moussouris \cite{Mou}, 
to quantum groups and the quantum link theory by Reshetikhin and Turaev \cite{RT},
Ohtsuki and Yamada \cite{OY}, Carter, Flath and Saito \cite{CFS} and others, 
to the theory of Lusztig's canonical bases 
\cite{Lu} by Khovanov and Frenkel \cite{FK}, Varchenko \cite{Var} and Frenkel, Varchenko and
Kirillov, Jr. \cite{FVK}.

Furlan, Stanev and Todorov \cite{FST} 
have extended outerplanar SL(2) invariants to the quantum 
algebra $U_q(sl(2))$ (for arbitrary spins). 
I became familiar with the description of the basis of the invariants of the tensor products 
of any irreducible representations of $SL(2)$ in the terms of the outerplanar graphs from 
Kuperberg's work \cite{K}. In \cite{M1}  
I gave a new proof of a classical theorem of 
Rumer, Teller and Weyl \cite{W} and its generalization for the case of arbitrary 
spins. In \cite{M2} I parametrized by outerplanar graphs the bases in the 
decompositions of any (repeated) tensor products of polynomial representations of 
$SL(2)$. Instead of the classical approach to the invariant theory using 
the straightening method, I used in \cite{M1,M2} the linear independence 
reason and the enumeration of the outerplanar graphs.

In recent work \cite{M3} I provided the spaces of invariants of tensor powers of 
the defining representations of $SL(n)$ 
with the bases parametrized by wave graphs introduced there 
especially for this purpose. The proof utilized a game similar to Tetris,  
named there L-tris, as well as the same linear independence reason as for the case
$n=2$ and the enumeration of wave graphs. 

Here I give similar constructions for 
$Sp(2n)$, parametrizing the invariants of tensor powers of the defining 
representation of $Sp(2n)$ by symplectic wave graphs introduced here especially    
for this purpose. The proof utilizes a game similar to L-tris,  
named here S-tris, as well as the same linear independence reason 
and the enumeration of symplectic wave graphs. 

I am preparing an article providing the space of invariants of tensor powers of 
the defining representations of orthogonal groups with the basis parametrized 
by (odd or even) orthogonal wave graphs introduced there especially for that 
purpose. Since we have a few different kinds of wave graphs, I propose to add to  
the name of wave graphs introduced in \cite{M3} 
the adjective `linear', i.\ e.\ refer to them as linear wave graphs
and use the term {\em wave graphs} for all of them: linear, symplectic and 
odd or even orthogonal wave graphs, as well as, exceptional wave graphs 
for the exceptional Lie groups. 

A symplectic $2n$-wave graph is a graph with the vertices $1,2, \dots, m$, 
each connected component of which is a path of length $\geq 1$ 
(i.\ e. it can't be a point), edges of which can be drawn 
in the book with $n$ pages, i.\ e.\ $n$ copies of the upper half-plane, glued 
along $\RR$, such that the first edge of each connected component, $\{i_1i_2\}$, 
is drawn on the first page; each edge $\{i_ji_{j+1}\}$ consequent to the    
edge $\{i_{j-1}i_j\}$ drawn on $k$-th page, is drawn either on $(k+1)$-th or 
$(k-1)$-th page, if they exist, I mean that the edge consequent to the edge drawn on 
the first page, must be drawn on the second page and the edge consequent to the edge 
drawn on the $n$-th page, must be drawn on $(n-1)$-th page; the last edge of the path, 
$\{i_l i_{l+1}\}$, supposed to be drawn on the first page; we suppose also that 
$i_1<i_2<\dots<i_{l+1}$ and edges of our symplectic wave graph don't intersect.

Symplectic $2$-graphs are exactly $1$-regular outerplanar graphs, or linear 
$2$-wave graphs which is not surprising because $Sp(2)=SL(2)$. Here are $4$ of a 
total number of $14$ of symplectic $4$-wave graphs with $6$ vertices: 

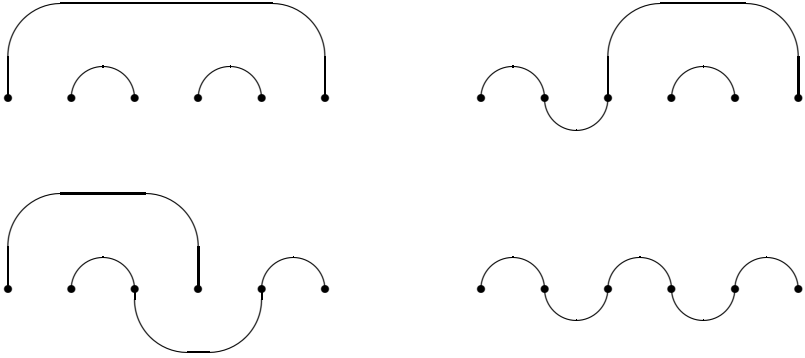
\begin{picture}(410, 72) 
  
\put(33,25){\circle*{3}}
\put(57,25){\circle*{3}}
\put(81,25){\circle*{3}}
\put(105,25){\circle*{3}}
\put(129,25){\circle*{3}}
\put(153,25){\circle*{3}}
\put(69,25){\oval(24,24)[t]}
\put(117,25){\oval(24,24)[t]}
\put(93,25){\oval(120,72)[t]}

\put(212,25){\circle*{3}}
\put(236,25){\circle*{3}}
\put(260,25){\circle*{3}}
\put(284,25){\circle*{3}}
\put(308,25){\circle*{3}}
\put(332,25){\circle*{3}}
\put(224,25){\oval(24,24)[t]}
\put(248,25){\oval(24,24)[b]}
\put(296,25){\oval(72,72)[t]}
\put(296,25){\oval(24,24)[t]}

\end{picture}

\begin{picture}(410, 62) 
  
\put(33,25){\circle*{3}}
\put(57,25){\circle*{3}}
\put(81,25){\circle*{3}}
\put(105,25){\circle*{3}}
\put(129,25){\circle*{3}}
\put(153,25){\circle*{3}}
\put(69,25){\oval(24,24)[t]}
\put(69,25){\oval(72,72)[t]}
\put(105,25){\oval(48,48)[b]}
\put(141,25){\oval(24,24)[t]}

\put(212,25){\circle*{3}}
\put(236,25){\circle*{3}}
\put(260,25){\circle*{3}}
\put(284,25){\circle*{3}}
\put(308,25){\circle*{3}}
\put(332,25){\circle*{3}}
\put(224,25){\oval(24,24)[t]}
\put(248,25){\oval(24,24)[b]}
\put(272,25){\oval(24,24)[t]}
\put(296,25){\oval(24,24)[b]}
\put(320,25){\oval(24,24)[t]}

\end{picture}

The corresponding invariants are 
\begin{gather}\label{1}
(\omega\otimes\omega\otimes\omega)^{(26543)},\quad 
((\omega\wedge\omega)\otimes\omega)^{(465)}\\
((\omega\wedge\omega)\otimes\omega)^{(1235)(46)},\quad 
\omega\wedge(\omega\otimes\omega) 
\label{2}\end{gather} 
where
\begin{equation}\label{3}
\omega=p_1\wedge q_1 + p_2\wedge q_2
\end{equation}
and for $\sigma\in S_6$, $t\in V^{\otimes 6}$ where $V$ is the defining 
representation of $Sp(4)$, the tensor $t^\sigma$ is the result of permutation 
$\sigma$ applied to the components of $t$; also, for $2$-tensor $\alpha$ and 
$(m-2)$-tensor $\beta$ we define $m$-tensor  
\begin{equation}\label{3.1}
\alpha\wedge\beta=\sum_{1\leq i<j\leq m} (-1)^{i+j-3}(\alpha\otimes\beta)^
{\sigma_{ij}}
\end{equation}
where $\sigma_{ij}\in S_m$ is the permutation mapping $1$ to $i$, $2$ to $j$ and 
other elements to the vacant places in increasing order, i.\ e.\ 
$\sigma_{ij}(3) < \sigma_{ij}(4) < \dots < \sigma_{ij}(m)$. 

For a symplectic wave graph $G$, denote $t_G$ the analogous tensor products of the 
basic invariants corresponding to the connected components, see Definition 
\ref{def2}. 

\begin{thm1}
Tensors $t_G$ parametrized by all $2n$-wave graphs with $m$ vertices, form a basis 
in the space of $Sp(2n)$-invariants in $V^{\otimes m}$, where $V$ is the 
$2n$-dimensional space of the defining representation of $Sp(2n)$.
\end{thm1}

The proof uses a game similar to Tetris, named here S-tris, linear independence 
reason, explicit formulas for the invariants and the enumeration of symplectic 
wave graphs. 

\section{The main theorem}
In this section we give all the necessary definitions and prove the main theorem.

Let $f$ be a field of characteristic $0$ and $Sp(2n)$ ---the group of 
$2n \times 2n$ $f$-matrices acting on $2n$-dimensional 
linear $f$-space $V$ with basis $B_{2n}=(p_1, \dots, p_n, q_1, \dots, q_n)$ 
by the standard way, preserving the symplectic 2-form 
\begin{equation}\label{4}
\omega=p_1\wedge q_1+\dots+p_n\wedge q_n .
\end{equation}

Recall some fundamental facts about the 
representations of $Sp(2n)$, see \cite{Hum}. The word {\em representation} will mean below 
a polynomial finite dimensional linear representation over $f$.
Every representation 
of $Sp(2n)$ is equivalent to a sum of irreducible representations. 
All classes of equivalence of the irreducible 
representations are parametrized by partitions of length $\leq n$. Denote 
$P_n$ the set of partitions of length $\leq n$ and
denote $\tilde{\rho}_\lambda$ the irreducible representation of $Sp(2n)$ 
corresponding to a partition $\lambda\in P_n$. 
Then $\tilde{\rho}_0$ is a trivial representation of dimension $1$ 
and $\tilde{\rho_1}$ is the standard representation in $V$ mentioned above.
To describe the decomposition of tensor products 
of irreducible representations, we'll use Young diagrams.

The same as in \cite{M3}, let us draw the Young diagrams rotated by $90^\circ$ 
counterclockwise.
Then we can interpret a Young diagram of a partition of length $\leq n$ as a 
Tetris position on a Tetris game field of width $n$, with non-increasing height 
of columns (from left to right). 

\begin{defn}\label{def1}
For a partition $\mu$ of length $\leq n$, 
denote $T_n(\mu)$ the set of partitions, Young diagram 
of which can be obtained from the Young diagram of $\mu$ by either 
dropping to it a $1\times 1$ block, or taking a top $1\times 1$ block in one 
of the columns of the Young diagram for $\mu$ and raising it back above 
the top of the Tetris game field. 
\end{defn}

Note that in contrast to the L-tris, defined in \cite{M3} for the description 
of the tensor products of representations of $SL(n)$, we don't delete 
complete Tetris rows here. 

Then
\begin{equation}\label{5}
\tilde{\rho}_\mu\otimes\tilde{\rho}_1\simeq\sum_{\lambda\in T_n(\mu)}\tilde{\rho}_\lambda .
\end{equation}

\begin{lemma}\label{lem1}
\begin{equation}\label{6}
\tilde{\rho_1}^{\otimes m}\simeq\sum_{\substack{\lambda\in L_n, |\lambda|\leq m\\
|\lambda|\equiv m \bmod 2}}\tilde{f}^\lambda_m(n) \tilde{\rho}_\lambda,
\end{equation}
where $|\lambda|$ denotes the weight (i.\ e. the sum of all parts) of 
a partition $\lambda$ and 
$\tilde{f}^\lambda_m(n)$  
is the number of symplectic lattice words in the alphabet $C_n=\{1, 2, \dots, n, 
\overline{1}, \overline{2},\dots, \overline{n}\}$ of length $m$ and weight 
$x_\lambda = x_1^{\lambda_1} \dots x_n^{\lambda_n}$, 
where a word $i_1\dots i_m$ is called a symplectic lattice word iff the weight  
of each its initial subword $x^{\sgn i_1}_{|i_1|}\dots x^{\sgn i_k}_{|i_k|}$ 
equals $x^\tau=x_1^{\tau_1}\dots x_n^{\tau_n}$ for a partition $\tau$, i.\ e.\ 
$\tau_1 \geq \dots \geq \tau_n \geq 0$ where $\sgn i=1,\thickspace
\sgn \overline{i}=-1,\thickspace |i|=i,\thickspace |\overline{i}|=i$ for 
$i\in C_n$.
\end{lemma}
\begin{proof} 
By induction on $m$, by iteration of \eqref{5}, 
the Young diagrams of the partitions $\lambda$ in the right hand side of \eqref{6} 
can be obtained by dropping or raising $m$ $1\times1$ blocks on the 
Tetris game field as were described above. Writing each time when a block drops or 
raises the number of the column where it drops, or overlined number of the 
column from the top of which it raises, one obtains a symplectic lattice word, 
because the definition of the symplectic lattice word means exactly that we have 
a Young diagram on each step of our game. 
\end{proof} 

\begin{cor}\label{cor1}
The dimension of the space of $Sp(2n)$-invariants in $V^{\otimes m}$ where $V$ is 
the defining representation of $Sp(2n)$, equals $\tilde{f}^0_m(n)$, 
the number of balanced 
symplectic lattice words in the alphabet $C_n$, where balanced means that 
the word contains the same number of 
$i$'s and $\overline{i}$'s for every $i$ from $1$ to $n$. 
\end{cor}
\begin{proof}Since the dimension of the space of invariants is a coffecient 
at $\tilde{\rho}_0$ in \eqref{6}, it equals $\tilde{f}^0_m(n)$ by Lemma \ref{lem1}.
\end{proof}
 
\begin{lemma}\label{lem2}
The subspace of $Sp(2n)$-invariants 
of $V\otimes V$ is one-dimensional and we can choose the fundamental form $\omega$ 
defined in \eqref{4} as a basis element in this subspace.
\end{lemma}
\begin{proof}
By the definition of $Sp(2n)$, $\omega$ is invariant.
By Corollary \ref{cor1}, the dimension of the space of invariants is equal to 
$\tilde{f}^0_2(n)=1$ since there is exactly one possible S-tetris game leaving nothing 
after two steps: drop the $1\times 1$ block at the first column and then raise it. 
$1\overline{1}$ is the corresponding unique balanced 
symplectic lattice word of length $2$. 
\end{proof}

\begin{defn}\label{def3}
A symplectic $2n$-wave graph is a graph with the vertices $1,2, \dots, m$, 
each connected component of which is a path of length $\geq 1$ (i.\ e. it can't be 
a point), edges of which can be drawn 
in the book with $n$ pages, i.\ e.\ $n$ copies of the upper half-plane, glued 
along $\RR$, such that the first edge of each connected component, $\{i_1i_2\}$, 
is drawn on the first page; each edge $\{i_ji_{j+1}\}$ consequent to the    
edge $\{i_{j-1}i_j\}$ drawn on $k$-th page, is drawn either on $(k+1)$-th or 
$(k-1)$-th page, if they exist, I mean that the edge consequent to the edge drawn on 
the first page, must be drawn on the second page and the edge consequent to the edge 
drawn on the $n$-th page, must be drawn on $(n-1)$-th page; the last edge of the path, 
$\{i_l i_{l+1}\}$, supposed to be drawn on the first page; we suppose also that 
$i_1<i_2<\dots<i_{l+1}$ and edges of our symplectic wave graph don't intersect.
\end{defn}

\begin{lemma}\label{lem3}
The number of symplectic $2n$-wave graphs with $m$ vertices is nonzero iff 
$m$ is even in which case it equals $\tilde{f}^0_m(n)$.
\end{lemma}
\begin{proof}
We'll construct a bijection between the set of balanced symplectic lattice words 
of length $m$ in the alphabet $C_n$ and the set of symplectic $2n$-wave graphs with 
$m$ vertices. After that Lemma \ref{lem3} will follow from Corollary \ref{cor1}.  

First, construct a mapping from graphs to words. For each vertex $i$ of a 
symplectic $2n$-wave graph denote $\alpha_i=k\in C_n$ such that $|k|$ is 
the largest number of pages of the book containing edges ending in $i$ and 
$\sgn k=\sgn (j-i)$ where $\{ij\}$ is the edge lying on the $k$-th page; in other 
words, $\sgn k=1$ for the initial vertex of the path and if running along the 
path through $i$ we come from a page with a smaller number to a page with a larger 
number, or $\sgn k=-1$ for the last vertex of a path and if running along the 
path through $i$ we come from a page with a larger number to a page with a smaller 
number. Then the word $\alpha_1 \dots \alpha_m$ must be 
a balanced symplectic lattice 
word. Since the weight of an initial subword of this word is a product of weights 
of initial subwords of paths, it is enough to show that the word corresponding to 
a path is a symplectic lattice word. 

For paths, we'll use induction on their lengths. 
There is the unique path of length $1$, the corresponding word $1\overline{1}$ is 
a balanced symplectic lattice word. Suppose that the words corresponding to paths 
of length less than $l>1$, are balanced symplectic lattice words. Take a path of 
length $l$. 

If it doesn't have other edges on its first page except the first edge 
and the last one, then deleting the first and the last vertices and the first page, 
we obtain a $2(n-1)$-wave path of length $l-2$ (after the appropiate renaming of the 
pages and vertices). By supposition, the word $2\dots \overline{2}$ 
corresponding to it must be a balanced symplectic lattice word in the alphabet 
$C_n\setminus \{1, \overline{1}\}$ containing only one $2$ (initial) and only one 
$\overline{2}$ (final). Thus, the weight of each initial subword must be of 
type $x_2x_3\dots x_k$ for some $k\geq 2$, or $1$ if this subword is the whole word. 
Thus, after adding deleted vertices $1$ and $\overline{1}$, the weights of the 
initial subwords will be $x_1x_2x_3\dots x_k$ for some $k$ or $1$ for the whole 
word that means that it will be a balanced symplectic lattice word. 

If the path of length $l$ contains other edges than the first and the last on 
its first page, let one of them be the edge $\{i(i+1)\}$, then $\alpha_i=\overline{2},
\thickspace \alpha_{i+1}=2$. In this case the complete subgraph of the given path 
with vertices from $1$ to $i+1$, edges of which are drawn on the same pages,
is a symplectic $2n$-wave path and analogously for the complete subgraph 
with vertices from $i$ to $l+1$. The length of these paths is less than $l$, 
thus the words corresponding to them are balanced symplectic lattice words. 
For $j\leq i$, the initial subword $\alpha_1\dots\alpha_j$ of the original word, 
is the same as for the first of two new paths, thus the weight of it has the 
required form. In a new word $\alpha_{i+1}=\overline{1}$. 
Thus the weight of a subword $\alpha_1\dots\alpha_i$ is 
$x_1$, in both words, old and new, as well as the weight of the 
initial one-letter word $\alpha_i$ of another new word. Hence all the other 
weights of the initial subwords $\alpha_1\dots\alpha_j$ with $j>i$ of the 
original word, will be the same as for the initial subword $\alpha_i \dots 
\alpha_j$ of the second new word, thus they have the required form as well.
By induction, we have proven than the constructed above word is a balanced 
symplectic lattice word. 

Now we construct the inverse mapping, from words to graphs. 
Let $\alpha=\alpha_1 \dots \alpha_m$ be a balanced symplectic lattice word. 
To construct all the edges on the $k$-th page of our book, take the letters 
$\alpha_i$ of this word such that $|\alpha_i|=k$ or $|\alpha_i|=k+1$. 
Write them in the increasing order of their indices and rename all the letters 
$\alpha_i=k+1$ to $\overline{k}$ and all the letters $\alpha_i=\overline{k+1}$ to 
$k$. We get a letter in alphabet $\{k,\overline{k}\}$. The same as usual, 
the same as for outerplanar graphs \cite{M1, M2, M3, W}, draw the outerplanar 
graph on $k$-th page of our book with the chosen vertices: one way of 
doing that is to put a left bracket instead of $k$, right bracket instead of 
$\overline{k}$ and connect the corresponding left and right brackets.  
Doing that for all $k$ from $1$ to $m$, we obtain a graph which must be a 
symplectic $2n$-wave graph. Indeed, it is easy to check all the requirements. 
Also, by construction, these two mappings are mutually inverse. Thus, we 
constructed a bijection between the set of balanced symplectic lattice words 
of length $m$ in the alphabet $C_n$ and the set of symplectic $2n$-wave graphs with 
$m$ vertices.
\end{proof}

\begin{lemma}\label{lem3.3}
The number of connected symplectic $2n$-wave graphs with $m$ vertices
equals $c_{m-2}({\cal P}_n)$, the number of walks of length $(m-2)$ from $1$ to $1$ on 
the path ${\cal P}_{n-1}$ of length $(n-1)$, i.\ e.\ a simple graph with $n$ vertices 
$1, \dots, n$ and $(n-1)$ edges $\{i,i+1\}$ for $i$ running from $1$ to $(n-1)$.
\end{lemma}
\begin{proof}
The same as in the proof of Lemma \ref{lem3}, we'll construct a bijection. 
Let us think that the vertices of ${\cal P}_{n-1}$ are the numbers of the pages 
of a book in which our connected symplectic wave graph is drawn. The first edge of our 
graph is drawn on the first page, it means that in the beginning we are in the initial 
vertex $1$ of ${\cal P}_{n-1}$. The second edge is drawn on the second page: it 
corresponds to moving from $1$ to $2$ in ${\cal P}_{n-1}$; and so on, if the edge 
consequent to an edge drawn on $k$-th page, is drawn on $(k\pm 1)$-th page, we 
are moving from $k$ to $(k\pm 1)$ in ${\cal P}_{n-1}$. The last edge is drawn on the 
first page, it means that at the end of our walk on ${\cal P}_{n-1}$ we are 
returning to the vertex $1$. Conversely, for each walk $\alpha_1\dots\alpha_{m-1}$ 
from $1$ to $1$ on 
${\cal P}_{n-1}$, we can construct a connected symplectic $2n$-wave graph with $m$ 
vertices, drawing its $k$-th edge on $\alpha_k$-th page. These two mappings are mutually 
inverse. Thus we constructed a bijection.
\end{proof} 

By Lemma \ref{lem3}, 
we have the same number of symplectic wave graphs as we need. Let us 
construct the corresponding invariants. 

\begin{defn}\label{def2}
For a symplectic $2n$-wave graph $G$ having 2 vertices $\{1, 2\}$ and an edge 
between them, drawn on the first page,  
denote $t_G=\omega$. For a non-connected symplectic $2n$-wave graph $G=G_1 \coprod  
G_2$ define
\begin{equation}\label{d2eq1}  
t_G=(t_{G_1^o}\otimes t_{G_2^o})^\sigma
\end{equation}
where $G_1^o$ and $G_2^o$ are symplectic $2n$-wave graphs obtained from $G_1$ and 
$G_2$ by reindexing their vertices in the same order and $\sigma$ is the 
permutation putting the vertices of graphs $G_1^o$ and $G_2^o$ in their tensor 
product \eqref{d2eq1} on their correct positions in $G$. For a path $G$ with 
$m>2$ vertices 
with a balanced symplectic lattice word 
$1\beta\overline{1}$ where $\beta=\beta_1\dots\beta_{m-2}$ is a balanced symplectic 
lattice word in the 
alphabet $C_n\setminus \{1,\overline{1}\}$ define 
\begin{equation}\label{d2eq2} 
t_G=\omega\wedge t_B
\end{equation}
where $B$ is $2n$-wave graph corresponding to the word $(\beta_1-1)\dots(\beta_{m-2}-1)$ 
in the alphabet $C_n$ and for $2$-tensor $\alpha$ and 
$(m-2)$-tensor $\beta$ we define $m$-tensor  
\begin{equation}\label{3.11}
\alpha\wedge\beta=\sum_{1\leq i<j\leq m} (-1)^{i+j-3}(\alpha\otimes\beta)^
{\sigma_{ij}}
\end{equation}
where $\sigma_{ij}\in S_m$ is the permutation mapping $1$ to $i$, $2$ to $j$ and 
other elements to the vacant places in increasing order, i.\ e.\ 
$\sigma_{ij}(3) < \sigma_{ij}(4) < \dots < \sigma_{ij}(m)$. 
\end{defn}

More general,

\begin{defn}\label{dw}
For a $k$-tensor $\alpha$ and 
$(m-k)$-tensor $\beta$ we define $m$-tensor 
\begin{equation}\label{dw1}
\alpha\wedge\beta=\sum_{1\leq i_1<\dots<i_k\leq m} (-1)^{(i_1-1)+\dots+(i_k-k)}
(\alpha\otimes\beta)^
{\sigma_{i_1\dots i_k}}
\end{equation}
where $\sigma_{i_1\dots i_k}\in S_m$ is the permutation mapping $1$ to $i_1$, 
$2$ to $i_2$, $\dots$, $k$ to $i_k$  and 
other elements to the vacant places in increasing order, i.\ e.\ 
$\sigma_{i_1\dots i_k}(k+1) < \dots < \sigma_{i_1\dots i_k}(m)$. 
\end{defn}

\begin{lemma}\label{lem3.2}
The wedge products of tensors defined above is associative, distributive 
respective to the addition and satisfies
\begin{equation}\label{dw2}
\alpha\wedge\beta=(-1)^{k(m-k)}\beta\wedge\alpha
\end{equation}
for a $k$-tensor $\alpha$ and $(m-k)$-tensor $\beta$.
\end{lemma}
\begin{proof}
Associativity and distributivity follows directly from Definition \ref{dw}. 
For \eqref{dw2},
note that the sign of 
an item of the sum in \eqref{dw1} coincides with the sign of the corresponding 
permutation $\sigma_{i_1\dots i_k}$ since this permutation has exactly 
$(i_1-1)+\dots+(i_k-k)$ inversions. Permutations in the left hand side and 
the right hand side of \eqref{dw2} differs on $\sigma_{i_1\dots i_k}$ with 
$i_j=m-k+j$ having $k(m-k)$ inversions.
\end{proof}

\begin{lemma}\label{lem5}
For any symplectic $2n$-wave graph $G$, 
in the lexicographical order of the monomial basis of $V^{\otimes m}$ corresponding 
to the ordering $p_1 < p_2 < \dots < p_n < q_n < \dots < q_2 < q_1$ of $B_{2n}$, the monomial 
\begin{equation}\label{c2eq1}
b_{\alpha(G)}=b_{\alpha_1}\otimes\dots\otimes b_{\alpha_m}
\end{equation}
where $\alpha(G)$ is a balanced symplectic lattice word corresponding to the 
symplectic $2n$-wave graph $G$ and 
\begin{equation}\label{c2eq2}
b_i=\begin{cases} p_i & \text{if $\thickspace \sgn i=1$,}\\ q_{|i|} & 
\text{if $\thickspace \sgn i=-1$,}  
\end{cases}\end{equation}
is the minimal monomial with a non-zero coefficient in $t_G$.
\end{lemma}
\begin{proof}
The proof is not very simple and we'll do it in a few steps. 
First, recall that 
\begin{equation}\label{l5eq1}
\omega^{\wedge k}= 
\underbrace{\omega\wedge\dots\wedge\omega}_{k\thinspace \text{times}}= 
k!\sum_{\substack{1\leq i_1<\dots<i_k\leq n\\ \sigma\in S_{2k}}}
(-1)^{\inv \sigma}(p_{i_1}\otimes \dots\otimes p_{i_k}\otimes q_{i_k}\otimes\dots\otimes q_{i_1})^\sigma .
\end{equation}
For a connected symplected $2n$-wave graph corresponding to a walk $12\dots (k-1)k
(k-1)\dots 21$ according to the bijection constructed in Lemma \ref{lem3.3}, the 
invariant $t_G$ is $\omega\wedge\dots\wedge\omega$. 
It follows from \eqref{l5eq1} that the basis monomial corresponding to the word 
$12\dots k\overline{k}\dots \overline{2}\overline{1}$ is the minimal monomial with a 
nonzero coefficient for that case. 

By induction we can deduce that for other connected symplectic $2n$-wave graphs 
\begin{equation}\label{l5eq2}
t_G=\omega^{\wedge k}\wedge (t_{G_1}\otimes t_{G_2}) 
\end{equation}
for some $k$ and symplectic $2n$-wave graphs $G_1$ and $G_2$. 

Using \eqref{l5eq2}, we can prove by induction that if we write down all the letters 
$i$ and $\overline{i}$, in the same order, 
from the index word of a monomial with a non-zero 
coefficient in $t_G$, we get a word in the $2$-letter alphabet $\{i,\overline{i}\}$ 
such that the weight of each initial subword is either $x$ or $x^{-1}$; 
for any $i$ from $1$ to $n$. Note that in a word $\alpha_1\dots\alpha_{2l}$ 
in $2$-letter alphabet $\pm 1$ with this condition, the letters $\alpha_{2k-1}$ and 
$\alpha_{2k}$ have different signs for all $k$ from $1$ to $l$. 

Denote $M_m(2n)$ the set of balanced words of length $m$ in the alphabet $C_n$ 
satisfying the condition above, i.\ e.\ such that for every $i$ from $1$ to $m$ 
the word in $2$-letter alphabet $\{i,\overline{i}\}$ containing all the entries 
of $i$ and $\overline{i}$, in the same order, has the weight either $x$ or $x^{-1}$ 
for each initial subword. Note that if we have 
\begin{equation}\label{l5eq3}
\alpha_{2k-1}=\overline{i},\quad \alpha_{2k}=i
\end{equation}
in one of such words, then transposing these letters $\overline{i}$ and $i$ in the 
original word we obtain a word in $M_m(2n)$ less than original one in the 
lexicographical order of the words mentioned in the statement of the Lemma.  

Denote $M_m^+(2n)$ the subset of $M_m(2n)$ containing such words that for every 
$i$ from $1$ to $m$ the word in $2$-letter alphabet  $\{i,\overline{i}\}$
containing all the entries 
of $i$ and $\overline{i}$, in the same order, is $i\overline{i}\dots i\overline{i}$, 
i.\ e.\ with the first letter $i$ and alternating of the letters on every step. 

Define for every word $\alpha=(\alpha_1\dots\alpha_m)\in M_m^+(2n)$ its {\em pattern} 
as a word 
\begin{equation}\label{l5eq4}
\pat(\alpha)=\sgn\alpha_1\dots\sgn\alpha_m
\end{equation}
of length $m$ in the alphabet $\pm 1$. By definition of $M_m^+(2n)$, 
the sum of all entries of $\pat(\alpha)$ is $0$ and 
each initial subword of this word
has a nonnegative sum of its entries.

Conversely, for every word $\delta=\delta_1\dots\delta_m$ of an even length $m$ 
in the alphabeth $\pm 1$ with zero sum of entries, 
each initial subword of which 
has a nonnegative sum of entries, we define its lattice word 
\begin{equation}\label{l5eq5}
\lat(\delta)=\alpha_1\dots\alpha_m \in M_m^+(2n)
\end{equation}
assuming $\alpha_1=1$, $\alpha_m=-1$ and for other $k$ from $1$ to $m$   
\begin{equation}\label{l5eq6}
|\alpha_k|=\max \{s_k(\delta), s_{k-1}(\delta)\}, \quad \sgn \alpha_k=\delta_k
\end{equation}
where $s_k$ means the sum of the first $k$ entries.

By definition, the pattern of the word $\lat(\delta)$ is $\delta$. 
By induction on $k$, reading the word from the beginning to the end,
we can check that $\lat(\delta)$ is the smallest word in 
$M_m^+(2n)$ of the pattern $\delta$ 
with respect to the lexicographical order. Indeed, the smallest 
possible first letter is $1$. If $\delta_2=1$, we can't have the second 
letter $1$ again, because the entries of $1$ and $\overline{1}$ must alternate; 
thus the smallest second letter is $2$. Otherwise, if $\delta_2=-1$, the only 
possibility for the second letter is $\overline{1}$. Later, if we have $s_{k-1}=i$, 
$s_k=i+1$, we can't use $1, 2, \dots, i$ for the $k$-th letter, because 
they were used earlier once more than the corresponding overline numbers (each of 
them), thus the smallest possible choice is $(i+1)$. If we have $s_{k-1}=i$, 
$s_k=i-1$, we must use one of $\overline{1}, \dots, \overline{i}$ as the $k$-th 
letter and the smallest possible choice is $\overline{i}$. Thus,    
$\lat(\delta)$ is the smallest word in 
$M_m^+(2n)$ of the pattern $\delta$. 

Introduce the inverse 
lexicographical order on the patterns, i.\ e.\ lexicographical order corresponding 
to the ordering $1\prec -1$. It is easy to see that for the patterns 
$\epsilon \prec \delta$ we have $\lat (\epsilon) \prec \lat (\delta)$. Indeed, 
if $k$ is the smallest integer such that $\epsilon_k \prec \delta_k$, i.\ e.\ 
$\epsilon_k=1, \thickspace \delta_k=-1$, then by construction the words 
$\lat (\epsilon)$ and $\lat (\delta)$ have the same letters on the first $k-1$ 
places and 
\begin{equation}\label{l5eq6.1}
(\lat (\epsilon))_k < (\lat (\delta))_k 
\end{equation}
since 
\begin{equation}\label{l5eq7}
(\lat (\epsilon))_k \in \{1, 2, \dots, n\},\quad 
(\lat (\delta))_k \in \{\overline{1}, \overline{2}, \dots, \overline{n}\} .
\end{equation} 

Now look at the patterns of the indices of the monomials with non-zero coefficients 
in $t_G$ for a connected symplectic $2n$-wave graph $G$. For any pattern of the 
analogous monomials for the symplectic $2n$-wave graph $B$ defined in 
Definition \ref{def2}, the smallest pattern that it can give us for $t_G$ is not 
less than if we add $1$ in the beginning of it and $-1$ at the end. Indeed, 
if it has $i$ $1$'s at the beginning before the first $-1$, adding $1$ at the beginning 
gives $(i+1)$ $1$'s at the beginning; the same as adding $1$ before the first $-1$, 
and it gives a smaller pattern than one with $i$ $1$'s in the beginning obtained by 
adding $1$ after the first $-1$. Analogously, adding $-1$ before the last $1$ gives 
a larger word than the adding $-1$ at the end since the line of $-1$'s containing it 
becomes longer. The smallest pattern that we can obtain in this way is when we add 
$1$ at the beginning and $-1$ at the end to the smallest pattern for $B$. 

The smallest pattern for $B$ with the correponding lattice word $1\overline{1}$ is 
$1-1$ since $t_B=\omega$ in that case. We have $\lat (1-1)=1\overline{1}$.
Prove by induction on $m$ that the same is true in general, i.\ e.\ the smallest 
pattern of the basis monomials of $t_G$ with non-zero coefficients is 
$\pat (\alpha(G))$ where $\alpha(G)$ is the corresponding lattice word for a 
symplectic $2n$-wave graph $G$. Since for non-connected graphs we obtain the 
smallest pattern by combining the patterns of the connected components, it is enough 
to prove that for a connected $G$ supposing
by inductive hypothesis that the smallest pattern for $B$ is $\pat (\alpha(B)$. 
Since we can obtain $\pat (\alpha(G))$ from  $\pat (\alpha(B)$ by adding $1$ at the 
beginning and $-1$ at the end, this is true as we had shown in the previous 
paragraph. 

Now when we know that the smallest possible pattern is $\pat \alpha(G)$ and the 
smallest lattice word with this pattern is
\begin{equation}\label{l5eq8} 
\lat (\pat (\alpha(G))=\alpha(G) ,
\end{equation}
the only thing that we have to do on the last step of our proof of Lemma 
\ref{lem5}, is to check that the coefficient at $b_{\alpha(G)}$ is non-zero. 
By induction, we'll prove that this coefficient is positive. 

First, reading the word $\alpha(G)$ for a connected symplected $2n$-wave graph $G$, 
we see that all the odd digits $1, 3, \dots$ are located on odd places, all 
even digits $2, 4, \dots$ are on even places, all the overlined odd digits 
are on even places and all the overlined even digits are on odd places. 
Indeed, it is true for the first $1$, and the next letter after $k$ can be either 
$k+1$, or $k$; the next letter after $\overline{k}$ can be either $\overline{k-1}$ 
or $k$, i.\ e.\ the parity changes on every step according to our hypothesis. Thus, 
by induction, it is true.   
  
Using that, we can prove, again by induction, that for a connected $2n$-wave graph 
$G$ all the monomials from $M_m^+(2n)$ have non-negative coefficients. The same is 
true for $B$ as well since its connected components don't interlace. Thus all 
the items giving a monomial $b_{\alpha(G)}$ in the wedge product \eqref{d2eq2} 
have non-negative coefficients, it means that they can't cancel and at least one 
of them, obtained from $p_1\wedge q_1\wedge t_B$ by putting $p_1$ in the first 
place, $q_1$ in the last place and $b_{\alpha_+(B)}$ between them, where 
$\alpha_+(B)$ is a word obtained from $\alpha(B)$ by changing $k$ to $k+1$
and $\overline{k}$ to $\overline{k+1}$ for all $k$ from 1 to $(n-1)$,    
has a positive coefficient since the coefficient at $b_{\alpha(B)}$ in $t_B$ is 
non-zero by induction hypothesis, and  $b_{\alpha_+(B)}$ has the same coefficient 
because $t_B$ is $Sp(2n)$-invariant and the linear transformation 
\begin{equation}\label{l5eq9}
p_k\mapsto p_{k+1},\quad q_k\mapsto q_{k+1} 
\end{equation}
for $k$ from $1$ to $n-1$ and  
\begin{equation}\label{l5eq10}
p_n\mapsto p_1,\quad q_m\mapsto q_1 ,  
\end{equation}
is symplectic, therefore preserve $t_B$.  
\end{proof}

\begin{theorem}
Tensors $t_G$ parametrized by all $2n$-wave graphs with $m$ vertices, form a basis 
in the space of $Sp(2n)$-invariants in $V^{\otimes m}$.
\end{theorem}
\begin{proof}
Lemma \ref{lem3} and Corollary \ref{cor1} show that the number of 
symplectic $2n$-wave graphs with $m$ 
vertices is exactly the same as the dimension of the corresponding space of 
$Sp(2n)$-invariants. By Definition \ref{def2} and Lemma \ref{lem2}, since 
tensor product of invariants are invariant as well as the result of permutation 
of tensor factors, for any symplectic $2n$-wave graph $G$, tensor $t_G$ is 
$Sp(2n)$-invariant. 
Hence if we prove linear independence of the set of $t_G$, 
our theorem will be proven. 
The proof is completely analogous to the proof of the particular case $n=2$ and 
the analogous theorem for $SL(n)$ given 
in my articles \cite{M1,M2,M3}. 

Denote $B$ the standard basis of $V^{\otimes m}$, consisting of $(2n)^m$ tensor 
products $X_1\otimes\dots\otimes X_m$ with $X_1, \dots, X_m \in \{p_1, \dots, p_n,  
q_1, \dots, q_n\}$. Suppose that $B$ is ordered lexicographically corresponding 
to the ordering $p_1 < p_2 < \dots < p_n < q_n < \dots < q_2 < q_1$. 
By Lemma \ref{lem5}, $b_{\alpha(G)}\in B$ is the minimal element of $B$ with 
a non-zero coefficient in the decomposition of $t_G$. Note that for different 
graphs $G$ the lattice words $\alpha(G)$ are different. So, we have 
$\tilde{f}^0_m(n)$ elements $b_{\alpha(G)}$ ---one for each $G$. 

To prove the linear independence of the set of $t_G$, we can show that the rank 
of the $\tilde{f}^0_m(n)\times (2n)^m$ matrix of the coefficients of $t_G$
in the basis $B$, is equal to $\tilde{f}^0_m(n)$.
To do that, we can find a non-zero $\tilde{f}^0_m(n)\times \tilde{f}^0_m(n)$ 
minor of that matrix. 
Consider the $\tilde{f}^0_m(n)\times \tilde{f}^0_m(n)$ 
submatrix with rows numbered by $G$ ordered the same way 
as $b_{\alpha(G)}$ and columns corresponding to $b_{\alpha(G)}$. 
As we noticed above,  
$b_{\alpha(G)}$ is the first element with a nonzero coefficient in the row $G$ and  
this coefficient is non-zero by Lemma \ref{lem5}. 
So, this matrix is upper triangular with non-zero elements on the diagonal, 
therefore its determinant is not $0$, 
that completes the proof of the linear independence of $t_G$.
\end{proof}

\begin{acn}
I would like to thank my thesis formal and informal advisors, Fan Chung Graham 
and Alexandre Kirillov, the supervisor of the 
part of the research supported by ONR grant, Andre Scedrov, 
our Mathematics Department and Graduate Group Chairpersons, 
Dennis DeTurck and Chris Croke, as well as UPenn Professors Ching-Li Chai,
Ted Chinburg, Murray Gerstenhaber, Herman Gluck, Michael Larsen, 
David Shale, Stephen Shatz, Herb Wilf, and Wolfgang Ziller for
useful discussions, Greg Kuperberg and Ivan Todorov for useful references
and my Gorgeous and Brilliant Wife, Bette, for her total support and love.
\end{acn}

\end{document}